\documentclass[conference]{IEEEtran}
\usepackage{xcolor}
\usepackage{amsmath}
\usepackage{amsfonts}
\usepackage{amssymb}
\usepackage{algorithm}
\usepackage{algpseudocode}
\usepackage{graphicx} 
\usepackage{subcaption}

\makeatletter
\def\ps@IEEEtitlepagestyle{%
  \def\@oddfoot{\mycopyrightnotice}%
  \def\@evenfoot{}%
}
\def\mycopyrightnotice{%
   {\footnotesize 979-8-3503-7240-3/24/\$31.00 \copyright 2024 IEEE\hfill}
}

\ifCLASSINFOpdf
\else
\fi

\begin{document}

\title{Distribution Network Restoration: Resource Scheduling Considering Coupled Transportation-Power Networks}


\author{
    \IEEEauthorblockN{Harshal D. Kaushik\IEEEauthorrefmark{1}, \textit{Member, IEEE,} Roshni Anna Jacob\IEEEauthorrefmark{1}, \textit{Graduate Student Member, IEEE,} \\ Souma Chowdhury\IEEEauthorrefmark{2}, \textit{Senior Member, IEEE}, Jie Zhang\IEEEauthorrefmark{1}, \textit{Senior Member, IEEE,}}
    \IEEEauthorblockA{\IEEEauthorrefmark{1}\textit{The University of Texas at Dallas, Richardson, Texas, USA} \\
    \IEEEauthorrefmark{2}\textit{University at Buffalo, Buffalo,  New York, USA}\\
    \{harshal.kaushik, roshni.jacob, jiezhang\}@utdallas.edu, soumacho@buffalo.edu}
}

\maketitle

\begin{abstract}

Optimal decision-making is key to efficient allocation and scheduling of repair resources (e.g., crews) to service affected nodes of large power grid networks. Traditional manual restoration methods are inadequate for modern smart grids sprawling across vast territories, compounded by the unpredictable nature of damage and disruptions in power and transportation networks.
This paper develops a method that focuses on the restoration and repair efforts within power systems. We expand upon the methodology proposed in the literature and incorporate a real-world transportation network to enhance the realism and practicality of repair schedules. Our approach carefully devises a {reduced}
network that combines vulnerable components from the distribution network with the real transportation network. Key contributions include dynamically addressing a  coupled resource allocation and capacitated vehicle routing problem over a new reduced
network model, constructed by integrating the power grid with the transportation network. This is performed using network heuristics and graph theory to prioritize securing critical grid segments. A  case study is presented for the 8500 bus system.  



\end{abstract}


%
\IEEEpeerreviewmaketitle

\section{Introduction}
Grid restoration, critical in the aftermath of extreme events like storms and other disasters, demands a comprehensive approach to efficiently restore electrical power while minimizing the downtime. 
Traditional restoration methods reliant on manual intervention lack adequacy for modern smart grids sprawling across vast territories, compounded by the unpredictable nature and varying degrees of damage to components during crises, rendering conventional scheduling impractical. Moreover, disruptions in transportation networks further complicate matters, requiring a holistic solution that integrates both power grid and transportation network interactions. 

To confront this multifaceted issue effectively, a scheduling strategy operating within the framework of an actual transportation network, akin to the Capacitated Vehicle Routing Problem (CVRP), is imperative. The CVRP \cite{EscobarLinfatiToth2014, LimWang2005}, a classic problem in operations research, involves optimizing the routes and schedules of vehicles with limited capacities to serve a set of customers or locations, aiming to minimize total travel time. 
In our specific case, initially, resources, namely the repair crews, are assigned to the damaged nodes. Subsequently, the aim is to devise a schedule that minimizes both the overall travel time and repair duration within the context of power grid restoration. 
Integrating CVRP-based optimization techniques with power grid restoration planning aims to achieve a synergistic solution optimizing both transportation and repair processes, thereby reducing overall downtime and enhancing system resilience against extreme events.



Restoration in power networks encompasses both service and infrastructure restoration~\cite{force2022methods}. Service restoration involves adopting emergency control actions such as network reconfiguration and load shedding to minimize the impact of network failures and provide power to critical loads~\cite{RoshniSteveLiChowdhuryGelZhang2022, ZhangRoshniSteveChowdhury2023}. 
While service restoration improves the response of the power network to extreme events, it does not bring the network to its pre-disrupted state, which requires repairing damaged components. Addressing this aspect of restoration is crucial, considering the size of the network, unpredictability in the post-outage network state, and the logistics associated with resource allocation for repairs, in light of the transportation network constraints.

\begin{figure*}[t!]
    \centering
    \begin{subfigure}[b]{0.27\textwidth}
        \centering
        \includegraphics[width=\textwidth]{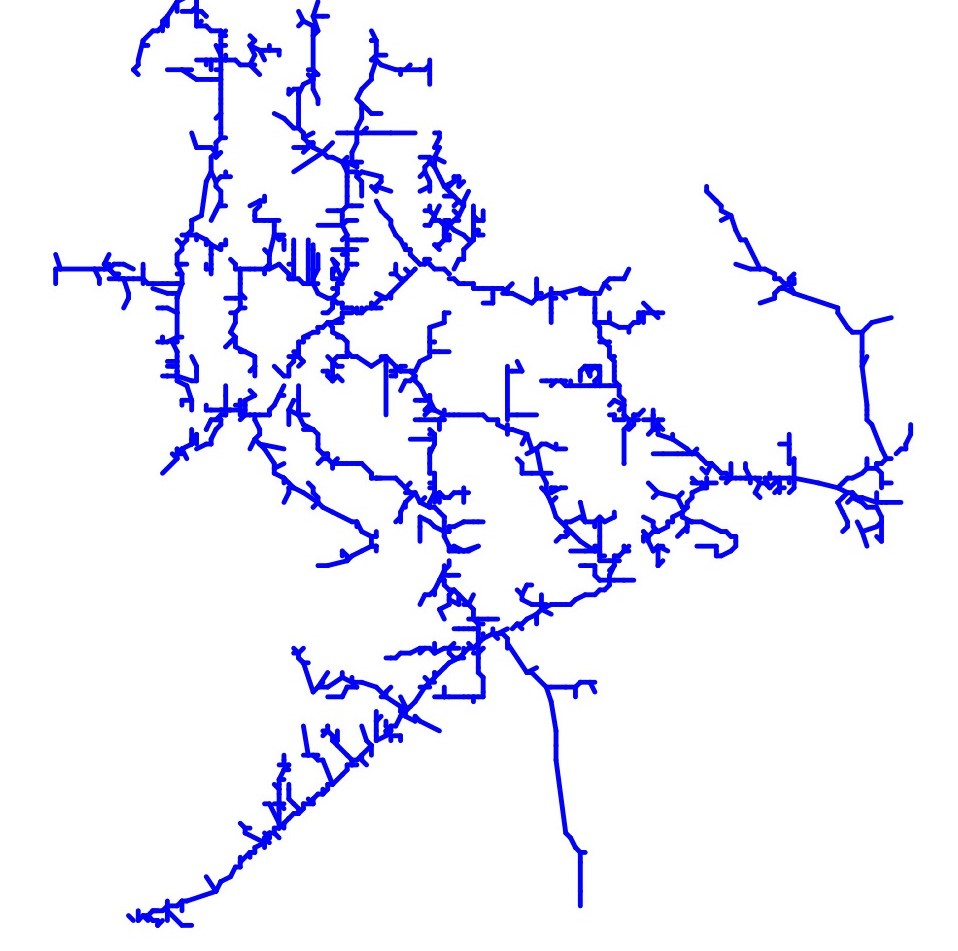}
        \caption{8500-node test feeder in OpenDSS}
        \label{fig:figure1-1}
    \end{subfigure}
    \hfill
    \begin{subfigure}[b]{0.345\textwidth}
        \centering
        \includegraphics[width=\textwidth]{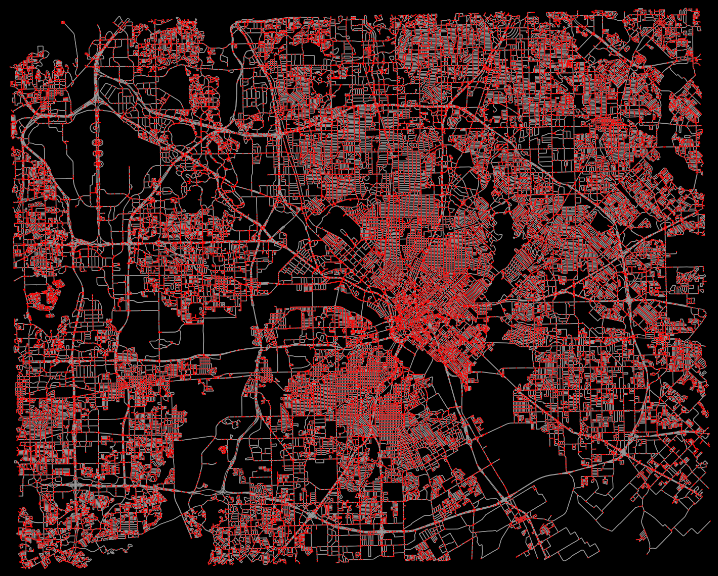}
        \caption{Transportation network of DFW area}
        \label{fig:figure1-2}
    \end{subfigure}
    \hfill
    \begin{subfigure}[b]{0.345\textwidth}
        \centering
        \includegraphics[width=\textwidth]{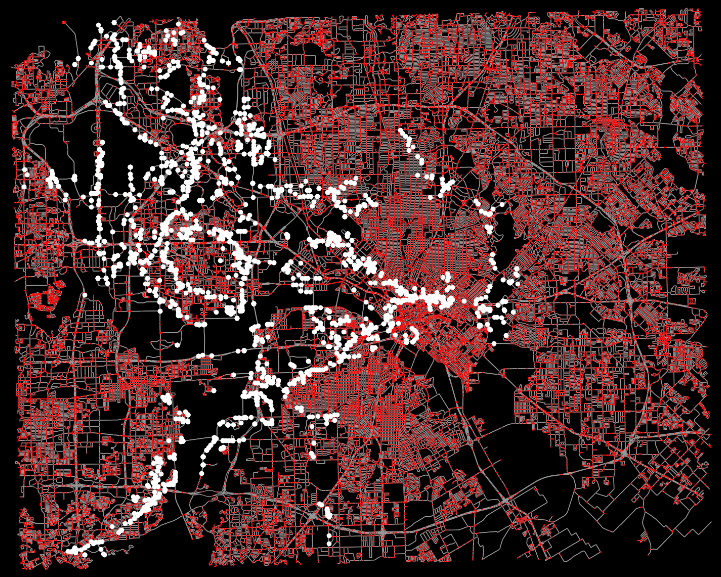}
        \caption{8500 nodes projected on transportation net}
        \label{fig:figure1-3}
    \end{subfigure}
    \caption{Projection of the IEEE 8500-node distribution test network on the DFW transportation network}
    \label{fig:figure1}
\end{figure*}
 Research efforts have explored the integration of repair and restoration within power transmission systems. Arab et al. \cite{ArabKhodaeiHanKhator2015} introduced pre- and post-hurricane models aimed at restoring power systems without the routing of crews. The routing of repair crews within transmission systems was included in \cite{HentenryckCoffrinBent2011}, where a multi-stage approach was utilized to separate routing and power-flow models. This methodology was further refined in \cite{SimonCoffrinHentenryck2012} through the application of a randomized adaptive vehicle decomposition technique. Zapata et al. \cite{ZapataSilvaGonzalezBurbanoHernandez2008} employed queuing theory and stochastic point processes to model repair schedules, dividing the distribution system into zones or service territories for individual resolution. Dynamic programming was considered in \cite{CarvalhoCarvalhoFerreira2007}, for directing the dispatch of repair crews and performing network reconfiguration. Researchers have also developed stochastic program \cite{XuGuikemaDavidsonNozick2007}, and large-scale optimization algorithms \cite{HDK2019, HDK2020, HDK2021} that are useful for crew assignment to inspect and repair damage, albeit overlooking crew routing considerations.

 The co-optimization of repairs and system restoration has been identified as a challenging problem \cite{HentenryckCoffrinBent2011}. The routing problem, a well-known NP-hard combinatorial optimization issue in operations research \cite{Laporte2009}, further compounds in complexity when integrated with emergent distribution system operation. Treating power system operation, restoration, and repair crew routing as independent problems is a common approach. However, practical interdependencies exist, such as the reliance on repairing damaged components for complete system restoration or leveraging distributed generators and automatic line switches to expedite restoration. Solely relying on utility operators' experiences for repair crew dispatch during outages may not yield optimal outage management plans. Consequently, there is a pressing need to design an integrated framework for optimally coordinating repairs and restoration efforts.

In \cite{ArifWangChenWang2}, the authors devised a mixed-integer linear program that considers the operation of the distribution network with crew routing. Subsequently, in \cite{ArifMaWangWangRyanChen2018}, their work was extended to address uncertainties via a stochastic mixed-integer linear program. They addressed this by decomposing the problem into two subproblems and employing parallel progressive hedging. Building upon their prior research in \cite{ArifWangChenWang2} and \cite{ArifMaWangWangRyanChen2018}, \cite{ArifWangChenWang1} enhanced their methodology by incorporating the three phase model considering unbalanced operation, modeling fault isolation constraints, and coordinating tree and line crews alongside resource logistics in the distribution system repair and restoration problem. This problem was solved over a complete graph network in parallel with the physical distribution network.

Our contribution is to extend this research of combining the restoration and repair. In this paper, we extend the methodology proposed in \cite{ArifWangChenWang2,ArifWangChenWang1} by considering a real transportation network. To make  more realistic and practical repair schedules, we have carefully devised a network that includes the most susceptible components from the distribution network and a real transportation network. {By borrowing the ideas from graph theory we have reduced down the nature of large-scale structure in a combined network.} Our contribution is summarized as follows:


\begin{itemize}
    \item 
    We introduce a network that intertwines the power grid and a real transportation network. Using graph theory, we model system connectivity during outages. Priority is given to securing the most susceptible and critical segments of the grid.
    \item We formulate the logistics and repair as a coupled resource allocation and CVRP  
    with constraints related to repair time, flow conservation, resource balancing, and repair crew attributes. This includes considering the repair time for each damaged component, as well as the travel time between damaged components for each crew. Additionally, crews are limited by the number of resources they can carry upon dispatch.
    \item We develop a two-stage iterative algorithm where in the first stage, we decide on the damaged nodes to be addressed first with the available resources {and calculate the weighted adjacency matrix for the selected nodes and depots}. Subsequently, in the second stage, we tackle the routing problem to minimize travel and repair time. This iterative approach focuses on eventually satisfying the demands at all the damaged nodes. 
\end{itemize}

\section{Model Formulation}\label{sec:model_formulation}

The CVRP is typically formulated over a graph comprising nodes and edges denoted as $G(V, E)$. We consider this undirected graph, the node set $V$ encompasses both the depot and the damaged components, while the edge set $E$, defined as $\{(i, j)|i, j \in V; i \neq j\}$. The graph $G$ is constructed by integrating a real transportation network and the power network. Next we discuss a stepwise approach to obtain $G$. 

\subsection{Integrated Transportation and Power Network}\label{subsec:2_A}
In this study, we employ the IEEE 8500-node distribution test feeder~\cite{IEEE8500} as a test case to construct an integrated network for our problem. 
We utilize the open source Distribution System Simulator (OpenDSS), to acquire the representation of the power network and the corresponding coordinates of the nodes. This test system is composed of both the primary and secondary levels, of which we consider the more critical primary level of the distribution network in our analysis. As a result we obtain a network comprising of 2,455 nodes and 2,454 edges, as shown in Fig. \ref{fig:figure1-1}. Subsequently, we translate the local coordinates of the primary distribution network into real geodesic coordinates.

For the transportation network, we opt to utilize the road network of the Dallas-Fort Worth (DFW) area. Leveraging Open Street Maps and NetworkX \cite{BoeingOSMNx}, we procure real road information, as shown in Fig. \ref{fig:figure1-2}. We then transform the nodes of the primary  level of the 8500-node feeder into the spatial framework of the transportation network. Next, we determine the nearest possible nodes by minimizing the Euclidean distances, {represented in Fig. \ref{fig:figure1-3}}. 
\subsection{Modified Resource Allocation and CVRP Formulations}
Our goal is to optimize the schedules and routes for each crew to reach the affected components, for the aim of maximizing the power restored and minimizing the total travel time. 
In this study, we employ a tailored blend of resource allocation and vehicle routing optimization models. As detailed in Algorithm \ref{alg:two_stage}, in the first stage we solve a resource allocation problem where we allocate appropriate repair crews to damaged nodes,  based on their requirements and priority. In the second stage, we determine the optimal routes for the allocated crews to minimize travel time.

Initially, the resource allocation problem assigns repair crews with the goal of minimizing repair time and maximizing power restoration, while adhering to vehicle capacity limits. Once crews are assigned, the second stage focuses on deriving optimal routes using the Capacitated Vehicle Routing Problem (CVRP) formulation. We will begin by delving into the specifics of the resource allocation approach next.

Resources here are the repair crews. These are allocated to the damaged nodes, such that maximal power is restored and the restoration times are comparatively smaller. {Power restored and repair time for node $i\in\mathcal{D}$ are denoted by $P_i$, and $T_i$, respectively.} Emphasis is primarily placed on maximizing power restoration, which is determined using OpenDSS powerflow simulations. The weighting of nodes based on the power capacity is considered by assigning suitable multiplication factors {$\alpha_1, \beta_1$} in the objective function in Equation \eqref{eq:1_RA} 
A continuous decision variable $y_{ik}$ denotes the total amount of resources allocated to the damaged node $i$ from crew $k$. Here, we denote the set of vertices $V \triangleq \{0, 1, \dots, N\},$ where $1, \dots, N-1$ are the damaged nodes and the start node and the end nodes are denoted by $0$ and $N$, respectively. Crew capacity is denoted as $Q_k$ for $k\in\mathcal{K}$.  Set $\mathcal{D}$ is the collection of all the damaged nodes. The {resource} demand at node $i\in\mathcal{D}$ is $q_i$. 
\begin{align}
\text{max}  & \sum_{k\in \mathcal{K}}\sum_{i\in {\mathcal{D}}} {\alpha_1}(P_{i}/ q_i) \ y_{ik} - {\beta_1}({T_{i }}/ q_i) \ y_{ik} \label{eq:1_RA} \\
\text{s.t.} \ & \sum_{i\in V }y_{ik} \leq Q_k &&  \hspace{-1.71cm} \text{ for } k \in \mathcal{K}   \label{eq:2_RA} \\
& y_{ik-1} \leq y_{ik} &&  \hspace{-1.8cm} \ \text{ for } \ i \in \mathcal{D}\label{eq:3_RA}\\
  & \ 0\leq y_{ik} \leq q_i &&  \hspace{-1.8cm}  \ \text{ for }  \ i \in \mathcal{D}, k\in \mathcal{K}.\label{eq:4_RA}
\end{align}


Constraint \eqref{eq:2_RA} ensures that the total capacity of the crew is not violated. Constraint \eqref{eq:4_RA} ensures the bound on the decision variable. Following the previous literature \cite{ArifWangChenWang1}, we add constraint \eqref{eq:3_RA} to maintain the crew sequence. For instance, there can be two primary crew types: line crews tasked with grid component repairs, and tree crews responsible for clearing obstacles at damage sites before the line crews begin repairs.

Subsequently, we delve into the CVRP formulation aiming to find optimal routes for repair crews to reach damaged components. Along with the travel time {(captured by matrix $c_{ijk}$ to travel from node $i$ to $j$ for crew $k$)}, our objective function in \eqref{eq:1} incorporates both the power restored {($P_{i}$)} and repair times {($T_{i}$)} with appropriate weights $\alpha_2$ and $\beta_2$. It is important to note that we focus on a single crew $k \in \mathcal{K}$ at any given time when determining these optimal routes, with 
$k$ serving as a parameter throughout the CVRP formulation.
Variable $x_{i,j,k}$ is a binary variable that indicates that crew $k$'s path traverses the edge $(i, j)$. This problem is solved over a graph $G(V,E)$, obtained from Subsection \ref{subsec:2_A}. 


\begin{align}
\min  & \sum_{i,j\in E} c_{ijk} x_{ijk} +\sum_{i\in\mathcal{D}}(\alpha_2 {T_{i }} -  \beta_2 P_i)\ t_i  \label{eq:1}\\ 
\text{s.t.} \ & \sum_{j\in V \text{ and }i,j\in E} x_{ijk} = 1  && \hspace{-0.75cm}i \in \mathcal{D} \ \text{ and }  \ q_{i,k}>0 \label{eq:2} \\
& \sum_{i\in V \text{ and }i,j\in E} x_{ijk} = 1  && \hspace{-0.75cm}j \in \mathcal{D} \ \text{ and }  \ q_{j,k}>0 \label{eq:3}  \\
 & \sum_{i\in V \backslash \{N\} } x_{ijk}\label{eq:6} \\ &\quad= \sum_{i\in V \backslash \{0\} } x_{jik}  && \hspace{-0.75cm}(i,j) \in E : j\in V \nonumber\\
& \text{if} \ x_{ijk}=1 \ \Rightarrow \ t_i + 1 = t_j &&\hspace{-0.75cm} (i,j) \in E\label{eq:7}\\
& \text{if} \ x_{ijk}=1 \ \Rightarrow \  t_j > t_i && \hspace{-0.75cm}(i,j) \in E\label{eq:8} \\
& \text{if} \ x_{ijk}=1 \ \Rightarrow \  u_i + q_{j,k} = u_j &&\hspace{-0.75cm} (i,j) \in E\label{eq:9} \\
& q_{i,k} \leq u_i \leq Q_k &&\hspace{-0.75cm} i \in V \label{eq:10} \\
 & x_{ijk} \in \{0,1\} &&\hspace{-0.75cm} (i,j) \in E\label{eq:11} \\
& u_{i} \in \mathbb{R} && \hspace{-0.75cm}i \in V \label{eq:12}\\
& t_{i} \in \mathbb{R} && \hspace{-0.75cm}i \in V. \label{eq:13}
\end{align}

\begin{figure*}[t!]
    \centering
    \begin{subfigure}[b]{0.297\textwidth}
        \centering
        \includegraphics[width=\textwidth]{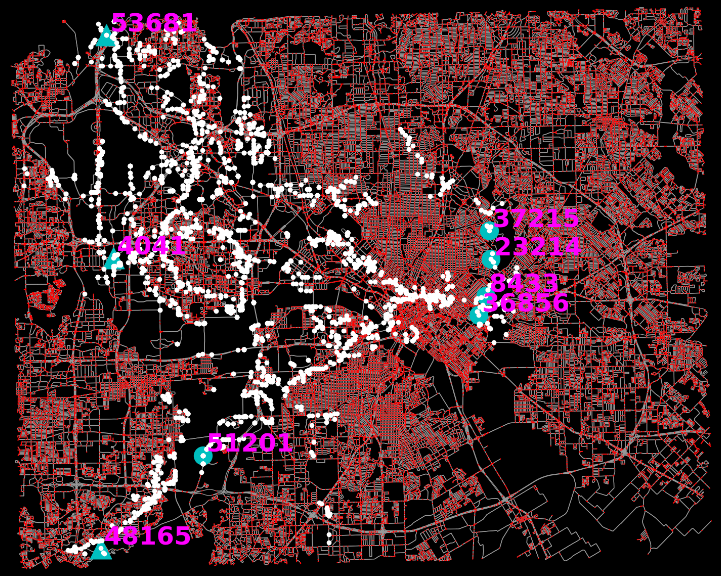}
        \caption{8500 net: damaged nodes and depots}
        \label{fig:figure2-1}
    \end{subfigure}
    \hfill
    \begin{subfigure}[b]{0.297\textwidth}
        \centering
        \includegraphics[width=\textwidth]{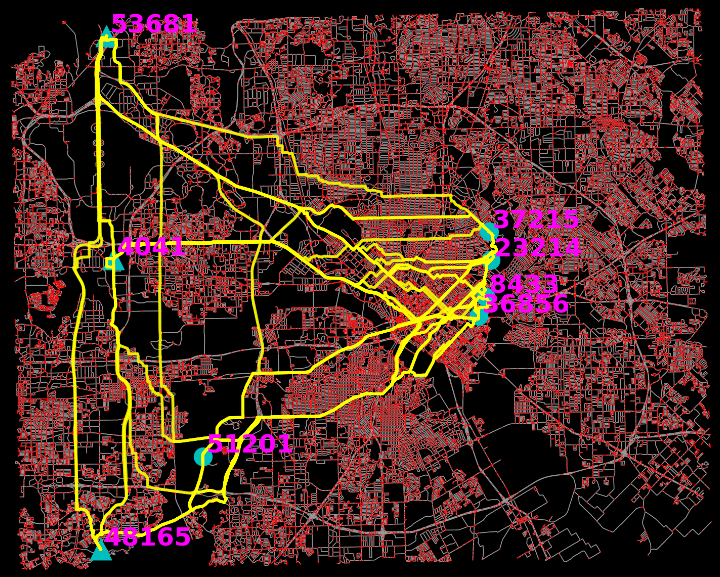}
        \caption{Shortest path among all the nodes}
        \label{fig:figure2-2}
    \end{subfigure}
    \hfill
    \begin{subfigure}[b]{0.333\textwidth}
        \centering
        \includegraphics[width=\textwidth]{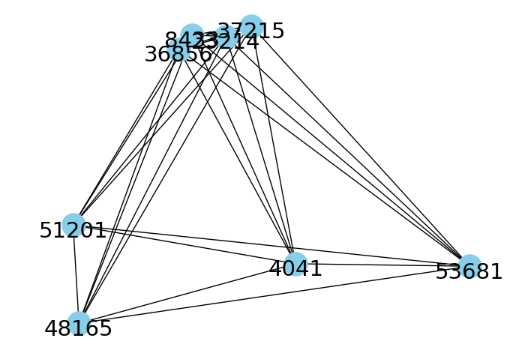}
        \caption{A complete graph from the damaged nodes}
        \label{fig:figure2-3}
    \end{subfigure}
    \caption{A complete graph constructed from five damaged nodes and three depots for the IEEE 8500-node test feeder overlaying the DFW transportation network.}
    \label{fig:figure2}
\end{figure*}

\begin{figure*}[t]
    \centering
    \begin{subfigure}[b]{0.447\textwidth}
        \centering
        \includegraphics[width=\textwidth]{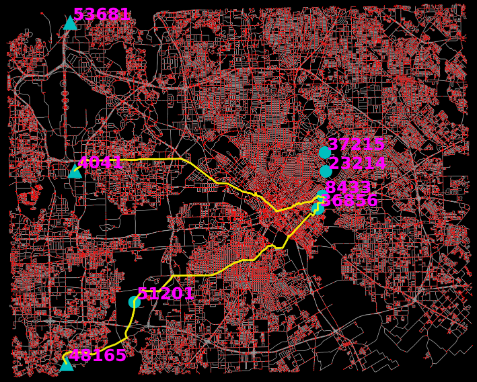}
        \caption{8500 net: damaged nodes and depots}
        \label{fig:figure3-1}
    \end{subfigure}
    \hfill
    \begin{subfigure}[b]{0.447\textwidth}
        \centering
        \includegraphics[width=\textwidth]{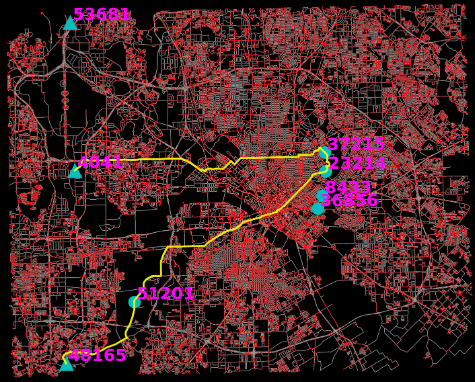}
        \caption{Shortest path among all the nodes}
        \label{fig:figure3-2}
    \end{subfigure}
    \hfill
    \caption{Solutions from the  iterations of Algorithm \ref{alg:two_stage} for the first experiment}
    \label{fig:figure3}
\end{figure*}

Constraints \eqref{eq:2} and \eqref{eq:3} enforce the requirement that each damaged node must be visited at least once during a run. 
Constraint \eqref{eq:6} represents the flow conservation constraint, dictating that once a crew arrives at a damaged component, it proceeds to the subsequent location upon completing repairs. Subtour elimination constraints \eqref{eq:7} and \eqref{eq:8} are incorporated to prevent the formation of suboptimal routes or ``subtours", thereby maintaining solution coherence and eliminating inefficient or redundant subroutes. Equation \eqref{eq:9} upholds the resource balancing constraint at each node $i \in V$. Additionally, constraints \eqref{eq:9} and \eqref{eq:10} ensure that the demand is met within the crew capacity.  Equations \eqref{eq:11}, \eqref{eq:12}, and \eqref{eq:13} make sure that $x$ is a binary variable, $u$ and $t$ are  in the real space.

Note that this CVRP is solved for every crew $k\in\mathcal{K}$. For every run, $k$ and $\mathcal{D}$ act as parameters. 
We explain the details of our solution procedure next in Algorithm \ref{alg:two_stage}.

\section{Algorithm Outline}
We present a two-stage algorithm to address our problem, tailored to different crews represented by $k \in \mathcal{K}$. These crews consist of two types: line crews,  and tree crews. The sequencing of crews is carefully managed to ensure optimal performance; specifically, tree crews are scheduled before repair crews to address obstacles efficiently. This allocation and proper sequencing of the resources, selection of damaged nodes $\mathcal{D}$, and the construction of the reduced adjacency matrix is done in the first stage. Second stage includes the solution of CVRP over an adjacency matrix of a complete graph connecting all the damaged nodes and depots. This process is repeated until all the demands at the damaged nodes $\mathcal{D}$ are satisfied. 
These steps are summarized in Algorithm \ref{alg:two_stage} below.


\begin{algorithm}
\caption{Optimal Vehicle Route Planning}\label{alg:two_stage}
\begin{algorithmic}[1]
\While{$q_i > 0$ for $i \in \mathcal{D}$}
\State {\bf Stage 1: Resource allocation and adjacency matrix}
\State  Decide the damaged node sequence to fulfil  $\mathcal{D}$.
\State Allocate crew resources with available capacity.
\State Obtain the weighted adjacency matrix for the reduced\quad complete graph.
\State {\bf Stage 2: Solve CVRP}
\State Solve CVRP and finalize the schedule for the crew.
\EndWhile
\end{algorithmic}
\end{algorithm}

\begin{figure*}[t]
    \centering
    \begin{subfigure}[b]{0.297\textwidth}
        \centering
        \includegraphics[width=\textwidth]{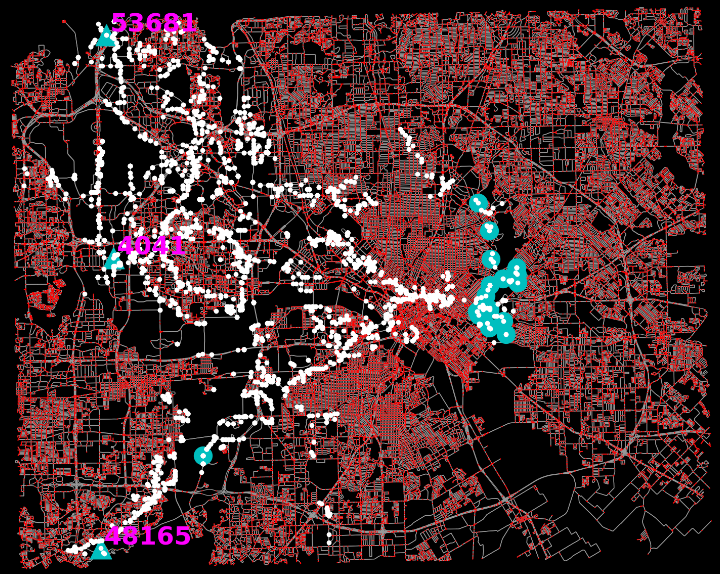}
        \caption{8500 bus: damaged nodes and depots}
        \label{fig:figure4-1}
    \end{subfigure}
    \hfill
    \begin{subfigure}[b]{0.297\textwidth}
        \centering
        \includegraphics[width=\textwidth]{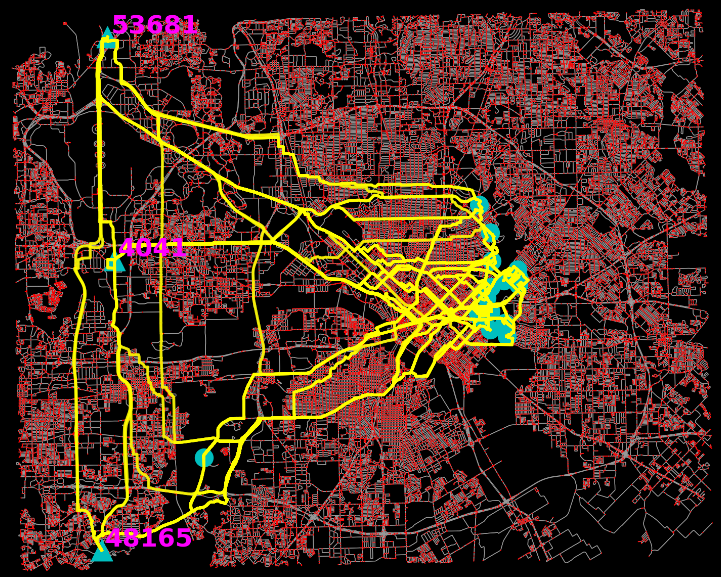}
        \caption{Shortest path among all the nodes}
        \label{fig:figure4-2}
    \end{subfigure}
    \hfill
    \begin{subfigure}[b]{0.333\textwidth}
        \centering
        \includegraphics[width=\textwidth]{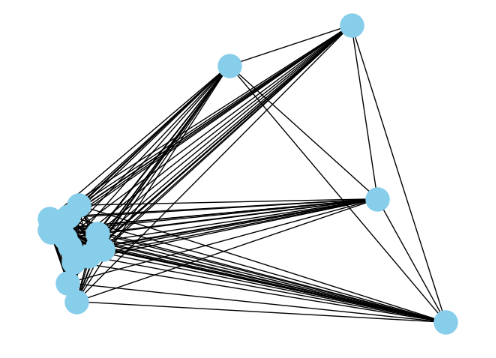}
        \caption{A complete graph from the damaged nodes}
        \label{fig:figure4-3}
    \end{subfigure}
    \caption{A complete graph constructed from 17 damaged nodes and three depots of 8500 bus, considering transportation network.}
    \label{fig:figure4}
\end{figure*}

\begin{figure*}[t!]
    \centering
    \begin{subfigure}[b]{0.447\textwidth}
        \centering
        \includegraphics[width=\textwidth]{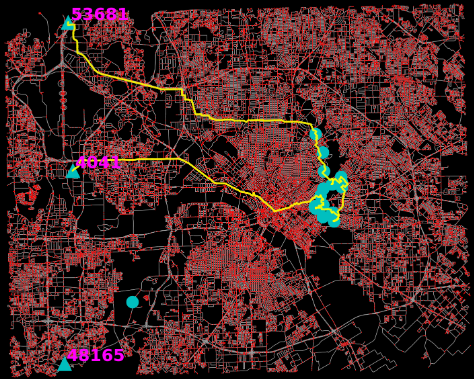}
        \caption{8500 bus: damaged nodes and depots}
        \label{fig:figure5-1}
    \end{subfigure}
    \hfill
    \begin{subfigure}[b]{0.447\textwidth}
        \centering
        \includegraphics[width=\textwidth]{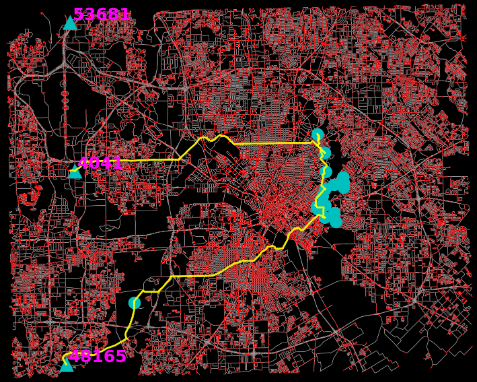}
        \caption{Shortest path among all the nodes}
        \label{fig:figure5-2}
    \end{subfigure}
    \hfill
    \caption{Solutions from the  iterations of Algorithm \ref{alg:two_stage} for the second experiment}
    \label{fig:figure5}
\end{figure*}


\section{Case Study: IEEE 8500-node test feeder on DFW transportation network}

In this section, we outline the solution process using a combined network of power and real roads around the DFW area. The process of construction of this network is described in Section \ref{sec:model_formulation}. For this case study,  we selected only the primary component of the 8500 bus network, which includes total  of 2455 nodes and 2454 edges. Initially, we identify the depots and damaged nodes, as depicted in Fig. \ref{fig:figure2-1}. Three depots and five damaged nodes have been considered in experiment 1. To create a complete graph, we compute the shortest path between all potential nodes and depots (Fig. \ref{fig:figure2-2}), followed by constructing the complete graph (Fig. \ref{fig:figure2-3}).  Table \ref{tab:power_restore_time} includes the power restoration potential provided the corresponding node is fixed (computed using OpenDSS), alongside the repair time (in hours) and the repair demand necessary to restore each damaged node.

\begin{table}[h]
\caption{Power restored and repair time of  damaged nodes}\label{tab:power_restore_time}
\centering
\begin{tabular}{|l|l|l|l|}
\hline
Node  & Power restored (kW) & Restore time (hours) & \begin{tabular}[c]{@{}l@{}}Repair\\ demand\end{tabular} \\ \hline
37215 & 78.43          & 3.59         & 6                                                       \\ \hline
23214 & 302.17         & 2.55         & 1                                                       \\ \hline
8433  & 10476.66       & 1.13         & 4                                                       \\ \hline
36856 & 10764.3        & 2.75         & 4                                                       \\ \hline
51201 & 10773.17       & 1.14         & 8                                                       \\ \hline
\end{tabular}
\end{table}

Following Algorithm \ref{alg:two_stage}, in the first stage, we computed all the weighted adjacency matrix entries for a complete graph. The crew's total capacity is set at 15 units. Initially, the repair prioritizes the last three nodes listed in Table \ref{tab:power_restore_time} because this combination yields the maximum power restoration, also taking into account corresponding repair time and repair demand. The output from the first iteration of Algorithm \ref{alg:two_stage} is shown in Fig. \ref{fig:figure3-1}. In the second iteration, the demand of the remaining nodes is fulfilled, and the outcomes are depicted in Fig. \ref{fig:figure3-2}. Algorithm \ref{alg:two_stage} halts once all demand is met.

In the second set of experiments, our goal was to evaluate the efficiency of Algorithm \ref{alg:two_stage} by including a significant number of damaged nodes, as depicted in Fig. \ref{fig:figure4-1}. A total of 17 nodes and 3 depots were chosen. The shortest paths for all 20 nodes were calculated, as depicted in Fig. \ref{fig:figure4-2}, and a complete graph was constructed, shown in Fig. \ref{fig:figure4-3}. For this particular run, we set the total crew capacity to 54. The repair demands of the damaged nodes are met within two iterations. The results of the first and second iterations are depicted in Fig. \ref{fig:figure5-1} and \ref{fig:figure5-2}, respectively. Once all the demands of the damaged nodes are fulfilled, Algorithm \ref{alg:two_stage} halts. The bottleneck lies in computing the shortest paths between all nodes and constructing the complete graph. Once the complete graph is obtained, solving the CVRP problem is not computationally challenging.

\section{Conclusion}
This paper focuses on the restoration and repair efforts within distribution grids. We incorporate a real transportation network to enhance the realism and practicality of repair schedules. Our approach carefully devises a reduced network that combines vulnerable components from the distribution network with the real transportation network, prioritizing critical grid segments. Notably, we have condensed a large transportation network to a manageable size suitable for conventional operations research solvers, allowing for effectively addressing a coupled resource allocation and capacitated vehicle routing problem formulation. A comprehensive analysis of algorithmic performance and implementation,  for different scenarios in the IEEE 8500-node distribution test feeder overlaying the DFW transportation network is performed. Implementation details, solution procedure, and the results are discussed in detail.

\section*{Acknowledgments}
This work was partially supported by the Office of Naval Research under ONR award number N00014-21-1-2530 and National Science Foundation under grant 2229417. The United States Government has a royalty-free license throughout the world in all copyrightable material contained herein. Any opinions, findings, and conclusions or recommendations expressed in this material are those of the author(s) and do not necessarily reflect the views of the Office of Naval Research and National Science Foundation.


\bibliography{KPEC.bib}

\end{document}